\documentclass[12pt,reqno]{amsart}

\usepackage[margin=1.15in]{geometry}
\usepackage{mathtools}
\usepackage{fontspec}
\usepackage{unicode-math}
\defaultfontfeatures{Ligatures=TeX,Scale=MatchLowercase}

\usepackage[final]{microtype}

\usepackage{amsthm}
\usepackage{thmtools}
\usepackage{bm}
\usepackage{mathrsfs}
\usepackage{amsmath}

\usepackage{booktabs}
\usepackage{enumitem}
\usepackage{tikz}
\usetikzlibrary{arrows.meta,positioning}
\usepackage{xcolor}

\definecolor{linkblue}{HTML}{145C9E}
\definecolor{citegreen}{HTML}{26734D}
\definecolor{urlmagenta}{HTML}{8A3FFC}
\usepackage[
  unicode,
  colorlinks=true,
  linkcolor=linkblue,
  citecolor=citegreen,
  urlcolor=urlmagenta
]{hyperref}

\numberwithin{equation}{section}

\newcommand{\R}{\mathbb{R}}

\newcommand{\tr}{\operatorname{tr}}

\newcommand{\End}{\operatorname{End}}

\title[]{A note on the BRST-BV Formalism Explanation of Deformation Quantization and Algebraic Index Theorem}

\author{Lihao Zhao}
\thanks{Artificial intelligence (AI) were used in the preparation of this manuscript to assist with mathematical computations and English-language editing.}
\address{School of mathematical sciences, East China Normal University, Shanghai, China}
\email{52300155001@stu.ecnu.edu.cn}

\date{November 2026}

\begin{document}

\begin{abstract}
In this note, we begin with considering the simple one-dimensional field theory describing a single free particle on manifold, and apply the BRST formalism and BV quantization to this model following the standard procedures in physics. Through calculations that, while not fully rigorous, have a more transparent physical interpretation, then subsequent geometric identifications, we recover the deformation quantization and algebraic index theorem on symplectic manifold. This note thus elucidates the physical origins and mechanisms of these mathematical results.
\end{abstract}

\maketitle

\tableofcontents

\section{Quick review of BRST and BV}

\subsection{Constrained system}

\par In a classical mechanical system, the action
$$
S[q]=\int_{t_0}^{t_1}L(q,\dot q,t)\,dt
$$
can be transformed into phase space via the Legendre transform:
$$
S[q,p]=\int_{t_0}^{t_1}(p_i\dot q^i-H(q,p,t))\,dt.
$$
Here, to invert normally to obtain $\dot q(q,p,t)$ and $H(q,p,t)$, the matrix
$W_{ij}=\frac{\partial^2L}{\partial\dot q^i\partial\dot q^j}$
must be invertible. Otherwise, such a system can only exist on certain submanifolds in phase space. These submanifolds $\Sigma$ are defined by the intersections of equations of the form
$T_a(q,p,t)=0$. However, for the time being, we arbitrarily extend $H$, which can only be inversely solved on the constraint surface, to the entire space.

\par This further introduces additional constraints. Since as a vector field $\{T_a,\quad \} \subseteq \ker(\omega|_\Sigma)$, the velocities along the directions of $\{T_a,\quad \}$ in the physical evolution on the submanifold may be arbitrary; this merely indicates that these directions are redundant. In other words, $\{T_a,\quad \}$ generates certain gauge symmetries. However, for such gauge symmetries to truly belong to the physical system defined on this constrained submanifold, the gauge transformation must ensure that points remaining on the constraint surface continue to lie on it—that is, it must always satisfy
$$
\{T,T_a\} = 0
$$
To meet this condition, a set of constraints $T$ is generally referred to as being first-class if locally
$$
\{T_a,T_b\} = f_{ab}^c(x) T_c.
$$
BRST theory aims to describe such systems using certain cohomology structures, whereas BV represents a scheme for quantizing such systems.

\par The challenge arises when $\{T,T_a\} = 0$ does not hold. A more favorable scenario is when $C_{\alpha\beta} = \{T_\alpha,T_\beta\}$ is invertible on the constraint surface, in which case it constitutes a set of second-class constraints. Under this condition, $\omega|_\Sigma$ ensures that the constraint surface itself becomes a symplectic manifold. There are two most common approaches to dealing with this situation. One involves replacing the Poisson bracket in full space with what is known as the Dirac bracket:
$$
\{F,G\}_D = \{F,G\} - \{F,T_\alpha\} C^{\alpha\beta} \{T_\beta,G\}
$$
which essentially corresponds to the Poisson bracket on the constraint surface.

\par Let the original constraints be denoted as $\theta_\alpha$. An alternative approach involves expanding the phase space by introducing a set of new variables $y^\alpha$ that satisfy
$$
\{y^\alpha,y^\beta\}=\Omega^{\alpha\beta},
\qquad \{x^i,y^\alpha\}=0,
$$
where $\Omega$ is a constant antisymmetric invertible matrix. The goal is to find a function
$$
T_\alpha(x,y)
=\theta_\alpha(x)+C_{\alpha\beta}(x)y^\beta+O(y^2)
$$
such that $\{T_\alpha,T_\beta\}=0$. The original function (including the Hamiltonian) $F(x)$ must also be extended to $\widetilde F(x,y)$, with the properties that
$$
\widetilde F(x,0)=F(x),
\qquad \{T_\alpha,\widetilde F\}=0.
$$
At this point, the system has been transformed into a set of first-class constraints.

\subsection{BRST cohomology}

\par \par Under the first-class constraints $T_a$, a set of new ghost variables $C^a$ and their conjugate momenta $P_a$ are introduced, satisfying
$$
C^aC^b=-C^bC^a,
\qquad (C^a)^2=0.
$$
$$
\{C^a,P_b\}
=\{P_b,C^a\}=\delta^a_b.
$$
We employ a grading scheme known as ghost number, defined as follows:
$$
\operatorname{gh}(x^i)=0,
\qquad \operatorname{gh}(C^a)=1,
\qquad \operatorname{gh}(P_a)=-1.
$$

\par The function space after continuous expansion is $C^\infty(M)\otimes\mathbb{C}[y,C,P]$. We construct the following function, namely the so-called BRST charge:
$$
Q=C^aT_a-\frac12C^aC^bf_{ab}^c  P_c
$$
Assume here $f_{ab}^c$ are constant. The BRST differential is defined as $s=\{Q,\quad\}$. It can be verified that since $\{Q,Q\}=0$, it follows that $s^2=0$. A brief computation of its actions yields:
$$
sf=C^a\{T_a,f\},
\qquad
sC^a=-\frac12f_{bc}{}^aC^bC^c,
\qquad
s P_b=T_b+C^af_{ab}{}^d P_d.
$$
If we consider $ C $ as gauge parameter, then $ s $ generates the gauge transformation. It is evident that the functions in $H^0(C^\infty(M)\otimes\mathbb{C}[y,C,P],s)$ are precisely those on the submanifold satisfying the constraints and gauge symmetry, i.e., the observables within this gauge-invariant constraint system. In other words, it is the space of gauge equivalence classes. The constraint equations vanish here because $s$ is exact, while the functions remain closed under $s$, thereby satisfying gauge symmetry. Further analysis shows that functions involving $P$ do not exist in $H^0$. This illustrates how BRST facilitates the description of such a system.

\par However, further work remains. If $f_{ab}^c$ is not constant, higher-order terms in $C,P$ must be introduced by termwise modification of the previously defined $Q$. Clearly, the new Hamiltonian $H_{\mathrm{BRST}}$ must satisfy
$$
\{Q,H_{\mathrm{BRST}}\}=0,
\qquad H_{\mathrm{BRST}}|_{C=P=0}=H.
$$
These conditions are also obtained by modifying $H$ termwise.

\subsection{Deformation quantization}
We now proceed to quantize this theory. Typically, a scheme known as deformation quantization is employed, under which the multiplication operation on the function space is transformed into a certain star product $\star$, and the Poisson bracket is transformed into the corresponding graded Lie bracket $[\quad,\quad]_\star$, satisfying the following relation:
$$
F\star G=FG+\frac{\hbar}{2} \{F,G\}+\hbar^2 C_2(F,G)+\cdots.
$$
This approach can be viewed as an extension of traditional canonical quantization, where
$$[q^i,p_j]=\hbar\delta^i_j$$
holds. Therefore, in Euclidean space, given the Poisson tensor $\pi^{ij}$ (the inverse of the symplectic form), the so-called Moyal product is commonly used:
$$F\star G=m\exp(\frac\hbar2\pi^{ij}\partial_i\otimes\partial_j)F\otimes G$$
Let $ m $ denote the replacement of the tensor product with the common multiplication. Additional variables are defined according to the relation:
$$[C^a,P_b]=\hbar\delta^a_b$$

\par The primary task involves systematically correcting $Q$ and $H_{BRST}$ term by term in powers of $\hbar$ within this quantum algebra, such that:
$$
Q_\hbar=Q+\hbar Q_1+\hbar^2Q_2+\cdots,
\qquad [Q_\hbar,Q_\hbar]_\star=0.
$$
$$
s_\hbar F=\frac1\hbar[Q_\hbar,F]_\star.
$$
Similarly, $H_\hbar$ is :
$$H_\hbar=H_{BRST}+\hbar H_1+\hbar^2H_2+\cdots,
\qquad [Q_\hbar, H_\hbar]_\star=0$$

\subsection{Gauge symmetry}
\par The quantization described above still fails to address the presence of gauge symmetry in the system, or equivalently, the issue of gauge transformations $\delta_\varepsilon =\varepsilon^a\{T_a,\quad \}$. Such symmetry effectively renders path integral quantization infeasible, as it forces the minima of the action to be degenerate. Consequently, the stationary-phase method cannot be employed in perturbative calculations of the path integral near these minima. If we wish to compute the partition function, deformation quantization does not seem to suffice, and we must instead resort to path integral quantization.

\par We now address gauge symmetry and re-perform the path integral quantization. First, we introduce a multiplier $u$ and rewrite the action as follows:
$$
S_{\mathrm{cl}}[q,p,u]
=p_i\dot q^i-H(q,p)-u^aT_a(q,p)
$$
On the constraint surface, this expression coincides with the original action. The key step is to replace the gauge parameter $\varepsilon^a$ with a ghost $C^a$; in this case, the infinitesimal gauge transformations coincide with the previously defined BRST differential $s$. To ensure gauge invariance of the action, we naturally impose the relation
$$su^a=-\dot C^a+f_{bc}{}^au^bC^c$$

\par The fundamental philosophy of BV is as follows: By introducing additional new variables known as antifields, suitable submanifolds are sought within the enlarged field space such that the path integral does not encounter degeneracy issues at extremum points on these submanifolds. Denote all previously defined variables by $\Phi^A$ and the field space by $\mathcal{F}$. For each variable, an antifield $\Phi^*_A$ is introduced, ensuring that their parity is opposite and the sum of their ghost numbers equals $-1$. This effectively expands the field space to $T^*[-1]\mathcal{F}$. The Poisson bracket is generalized to the BV bracket as follows:
$$
\{F,G\}_{\mathrm{BV}}
=\left(
\frac{\partial_{right}F}{\partial\Phi^A}
\frac{\partial_{left}G}{\partial\Phi_A^*}
-\frac{\partial_{right}F}{\partial\Phi_A^*}
\frac{\partial_{left}G}{\partial\Phi^A}
\right)
$$
More precisely, in field theory, this should be expressed as a functional derivative:
$$
\{F,G\}_{\mathrm{BV}}
=\int dt\left(
\frac{\delta_rF}{\delta\phi^a(t)}\frac{\delta_lG}{\delta\phi^*_a(t)}
-\frac{\delta_rF}{\delta\phi_a^*(t)}\frac{\delta_lG}{\delta\phi^a(t)}
\right).
$$
This graded bracket possesses a ghost number of 1.

\par Next, the action is expanded to an even function with ghost number 0, denoted $S$, such that
$
S|_{\Phi^*=0}=S_{\mathrm{cl}}
$. It is still required that the new $s_{\mathrm{BV}}=\{S,\quad\}_{\mathrm{BV}}$ satisfies $s_{\mathrm{BV}}^2=0$, which implies
${\frac12\{S,S\}_{\mathrm{BV}}=0}$. This equation is known as the classical master equation. In the simple cases considered here, one may directly take
$$
S=S_{\mathrm{cl}}-\sum_A(s_{BV}\Phi^A)\Phi_A^*
$$

\subsection{Gauge fixing}
\par In fact, in the antifield method, we no longer require the original $P$: instead, we use $-u^*$ to replace it. It can be verified that
$$
s_{\mathrm{BV}}(-u_b^*)
= T_b + c^a f_{ab}^d (-u_d^*).
$$
Thus, the variables in our field space include $q$, $p$, $u$, $C$, and their antifields. It is straightforward to observe that the term in $s_{\mathrm{BV}}$ that does not involve antifields corresponds to the action of $s$. Next, we impose a set of gauge-fixing conditions $\chi(\phi) = \chi(q, p, u) = 0$. The submanifold determined by these conditions should intersect every gauge orbit on the original constraint surface exactly once. Of course, we first assume the existence of such a selection. A simple approach is to perform a path integral over the intersection of this submanifold and the original constraint submanifold, thereby achieving gauge reduction. Since gauge transformations are already generated by $s$, this requirement is equivalent to demanding that the matrix $M^a_b$ in
$$
s\chi^a = M^a_b C^b
$$
be invertible.

\par The actual implementation in BV proceeds as follows. We introduce new variables: the odd anti-ghost $\bar{C}_a$, the even auxiliary field $b_a$, and their respective antifields. The following relations are imposed:
$$
\operatorname{gh}(\bar{C}_a) = -1, \quad \operatorname{gh}(b_a) = 0,
$$
$$
s\bar{C}_a = b_a, \qquad sb_a = 0.
$$
Furthermore, the action is extended to
$$
S_{\mathrm{BV}} = S_{\mathrm{cl}} - \sum_A (s_{\mathrm{BV}}\Phi^A)\Phi_A^* - b_a\bar{C}^{*a}.
$$
In the cohomology of classical differential $s_{\mathrm{BV}}=\{S_{\mathrm{BV}},\quad\}_{\mathrm{BV}}$, both pairs of variables and their corresponding antifields are contractible. As in previous studies, the space of physical observables seems to be zeroth cohomology space of $s_{\mathrm{BV}}$, though now it's somehow just classical observable. We will see that's not really the space of physical observables. 

\par Consider the so-called gauge-fixing fermion:
$$
\Psi=-\bar C_a\left(\chi^a(\phi)+\frac\xi2b^a\right)
$$
where $\xi$ is an arbitrary constant parameter. Here, we must define a matrix $h_{ab}$ such that $b_a=h_{ab}b^b$. The gauge fixing is implemented by restricting the field space to one of its Lagrangian submanifolds:
$$
\mathcal L_\Psi=
\left\{\Phi_A^*=\frac{\partial_l\Psi}{\partial\Phi^A}\right\}
\subset T^*[-1]\mathcal F
$$
This implies that a symplectic structure of degree $-1$ is assigned to the field space:
\[
\omega_{\mathrm{BV}}
=\sum_A\delta\Phi_A^*\wedge\delta\Phi^A.
\]
The corresponding Poisson bracket is then $\{\quad,\quad\}_{\mathrm{BV}}$. This procedure eliminates the antifield. It can be verified that the action after restriction, or gauge fixing, is given by:
$$
\begin{aligned}
  S_{\mathrm{gf}}&=S_{\mathrm{BV}}|_{\mathcal L_\Psi}=S_{\mathrm{cl}}-s\Psi 
  \\
  &=p_i\dot q^i-H(q,p)-u^aT_a(q,p)+b_a\chi^a+\frac\xi2b_ab^a -\bar C_aM^a{}_bC^b
\end{aligned}
$$

\par Performing a path integral over this action:
$$\int D[\phi]e^{\frac{i}{\hbar}S_\mathrm{gf}}$$
The Berezin integral of the terms involving $\bar C$ and $C$ yields $\det(M)$, known as the Faddeev–Popov determinant. This determinant represents the perturbation induced by gauge fixing on the integration volume. The term involving $b$, after completing the square, leads to a Gaussian integral giving $\exp(\frac{-\chi^2}{2\xi})$; taking $\xi \to 0$ corresponds to restricting to the gauge-fixed submanifold. The term involving $u$ also yields a Gaussian integral equal to $\delta(T)$, thereby restricting to the constraint submanifold (so that an arbitrary initial extension of $H$ is no longer an issue). What remains is the original expression $\int D[\phi]\exp(\frac{i}{\hbar}\int dt(p\dot q-H))$. In this manner, we have completed the path integral quantization of the original classical mechanical system using the BV formalism.

\subsection{Gauge independence}

\par However, we have not yet examined whether the BV integral remains invariant under these arbitrarily chosen gauge fixings. In other words, we have not determined whether the integral of $S_{\mathrm{BV}}$ over different Lagrangian manifolds $\mathcal L_\Psi$ yields the same value. Clearly, $\{\Phi^*_A=0\}$ represents the initial field space without antifield extension and gauge fixing; this implies that as long as the gauge remains unchanged, the BV integral we have computed is indeed the one originally sought. We argue that each Lagrangian manifold $\mathcal L_\Psi$ obtained from gauge-fixing fermions in the manner described earlier can be smoothly deformed from $\{\Phi^*_A=0\}$.

\par The subsequent discussion requires the so-called BV Laplacian:
$$
\Delta=
\sum_A(-1)^{|\Phi^A|}
\frac{\partial_l}{\partial\Phi^A}
\frac{\partial_l}{\partial\Phi_A^*}
$$
or in field theory, 
$$
\Delta=\int dt
\sum_A(-1)^{|\Phi^A|}
\frac{\delta_l}{\delta\Phi^A(t)}
\frac{\delta_l}{\delta\Phi_A^*(t)}
$$
It should be noted that, just as we need to specify a measure or volume form when defining the Laplacian in conventional contexts, the definition of this BV Laplacian already assumes that we have chosen the volume forms $\rho=\delta\Phi^A \delta\Phi^*_A$, which are equivalent to the odd symplectic forms discussed earlier. Under these assumptions, we have:
$$
\Delta(FG)
=(\Delta F)G+(-1)^{|F|}F\Delta G
+(-1)^{|F|}\{F,G\}_{\mathrm{BV}}.
$$
Crucially, for any $\mathcal L$, there exists a Stokes’ theorem-like formula:
$$
\int_{\mathcal L}\Delta f \sqrt{\rho}=0
$$
This is because it is always possible to select local Darboux coordinates such that $\mathcal L$ is locally $\{\Phi^*_A=0\}$. Now consider a family of gauge-fixing fermions $\Psi_t$, which define a smooth family of Lagrangian manifolds $\mathcal L_{\Psi_t}$. It can be shown that
$$
\{\dot \Psi_t,\Phi^A\}_{BV}=0,
\qquad
\{\dot \Psi_t,\Phi_A^*\}_{BV}
=\frac{\partial_l\dot\Psi_t}{\partial\Phi^A},
$$
meaning that $\dot \Psi_t$, acting as a Hamiltonian function, generates the flow described by its Hamiltonian vector field $X_{\dot \Psi_t}=\{\dot \Psi_t,\quad\}$ precisely through the deformation of $\mathcal L_{\Psi_t}$ over time. After further calculations, we obtain:
$$
\begin{aligned}
\frac{d}{dt}\int_{\mathcal L_{\Psi_t}}f\sqrt{\rho}
&=\int_{\mathcal L_{\Psi_t}}L_{X_{\dot \Psi_t}}f\sqrt{\rho} \\
&=-\int_{\mathcal L_{\Psi_t}}\Delta(\dot \Psi_t f)\sqrt{\rho}-\int_{\mathcal L_{\Psi_t}}\dot \Psi_t \Delta f \sqrt{\rho}
\end{aligned}
$$
Thus, it follows that as long as $\Delta(e^{\frac{i}{\hbar}S_{BV}})=0$, the BV path integral of $S_{BV}$ is independent of the gauge choice. This condition is known as the quantum master equation.

\par However, $S_{BV}$ generally fails to satisfy this condition. Therefore, we must modify $\hbar$ term by term as follows:
$$
S = S_{BV} + \hbar S_1 + \hbar^2 S_2 + \cdots
$$
so that $S$ satisfies the quantum master equation. In fact, there is a potentially more fundamental perspective: we are modifying the measure $\rho$ in field space to $e^{\hbar S_1 + \hbar^2 S_2 + \cdots}\rho$, thereby ensuring that the path integral becomes independent of gauge choices under this measure.

\subsection{Observables}

\par This also establishes the requirements for observable quantities. Consider the operator called quantum differential 
$$
\mathcal Q f=\{S,F\}_{BV}-i\hbar\Delta f
$$
Under the quantum master equation, it holds that
$$\Delta(fe^{\frac{i}{\hbar}S})=\frac{i}{\hbar}Q(f) e^{\frac{i}{\hbar}S}$$
In other words, only the integrals of observable quantities $f$ for which $Q$ closed remain independent of gauge choices, whereas the integral of such observable quantities vanishes when $Q$  exact. Analogous to the BRST structure discussed earlier, here $Q$ generates the gauge symmetries within the quantum system, and its zeroth cohomology corresponds to the space of physical observables.

\par In the discussions above, we sometimes conflated the notation and terminology used in finite-dimensional examples with those in infinite-dimensional scenarios (i.e., true field theory). For instance, the path integral referred to in finite dimensions as $\int dqdp \, e^{\frac{i}{\hbar}p\dot q}$ is actually a Gaussian integral in Euclidean space. The genuine path integral is given by $$\int D[q(t)]D[p(t)] \, e^{\frac{i}{\hbar}\int dt \, (p(t)\dot q(t))}.$$

\par In the following discussion, we assume that the reader is familiar with Fedosov's deformation quantization, and the algebraic index theorem in physics style. (ref.\cite{F96}\cite{FFS05}\cite{GLL17}). The discussion of the BRST formalism and deformation quantization in this paper partly overlaps with that in \cite{GL01}. Our treatment of BV quantization and the computation of the partition function is essentially equivalent to that of \cite{GLL17}\cite{GLX21}, although their approach is presented in a different manner. 

\par For the subsequent discussion, we make the substitution $i/\hbar\to 1/\hbar$ and reverse the sign of actions relative to the convention adopted in this chapter.

\section{Deformation quantization in BRST formalism}

\subsection{Motivation of the model}

\par We aim to demonstrate that for the simplest and most natural field theory, applying path integral quantization via BRST and BV theories in a systematic manner naturally leads to deformation quantization and the algebraic index theorem. Consider the action of a free particle in $\mathbb{R}^n$:
$$S[q(t)]=\frac12\int|\dot q|^2 dt.$$
By performing a Legendre transformation, we obtain the action in phase space $(T^*\mathbb{R}^{n},\omega=dq^i\wedge dp_i)$ as
$$S[q(t),p(t)]=\int p_i\dot q^i-\frac{1}{2}|p|^2 dt.$$
We may rescale the volume of $\mathbb{R}^n$, as theoretically, such a rescaling does not alter any physical phenomena. Specifically, if we transform the metric on the entire $\mathbb{R}^{n}$ by $g\to g/\ell$, then the metric along the fibers of $T^*\mathbb{R}^{n}$ must undergo an analogous global transformation, namely $g^*\to \ell g^*$. Taking the infinite-volume limit as $\ell\to 0$, the phase space action reduces to
$$S=\int p_i\dot q^i dt.$$
This expression can also be written in a manner independent of symplectic coordinates. By introducing a symplectic potential $\lambda=p_idq^i$ with $d\omega=\lambda$, and merging the coordinates on $T^*\mathbb{R}^{n}$ into $X=(q,p)$, we obtain
$$S=\int \lambda(X)_i\dot X^i dt.$$
This constitutes the model we will consider subsequently. This is equivalent to setting the Hamiltonian $ H = 0 $, and consequently, its equation of motion follows $\dot X=0$. Therefore, when subsequently performing path integral quantization via perturbation theory, we consider formal neighborhoods of constant loops.

\par Moreover, upon performing a Legendre transformation on the Lagrangian $\mathcal L(X,\dot X)=\lambda_i(X)\dot X^i$, it is observed that
$$
P_i
=
\frac{\partial \mathcal L}{\partial\dot X^i}
=
\lambda_i(X)
$$
contains no explicit dependence on the velocities. Consequently, it is impossible to invert this relation to obtain $\dot X(X,P)$; instead, a set of constraints $P_i-\lambda_i(X)=0$ is generated. To address these constraints, we must employ the BRST and BV formalism.

\subsection{One-dimensional free field theory on symplectic manifolds}
\par The one-dimensional field theory we consider is defined on a $2n$-dimensional symplectic manifold $M,\omega$. Additional information may also be incorporated: a rank-$r$ complex vector bundle $E$. The manner in which “particles move on the manifold” must be specified in advance, namely by fixing a torsion-free symplectic connection $\nabla^{TM}$ and a connection $\nabla^{E}$ on $E$. As we shall see later, it suffices to study fields within neighborhoods of classical trajectories—specifically, loops within neighborhoods of constant loops, or those loops that are sufficiently small. Therefore, we proceed directly by analyzing the situation in local coordinates $X^i$ within a coordinate chart $U$, where $\omega = d\lambda$ can be taken. The free action is given by:
$$
S
=-\int
\left(
\lambda_i(X)\dot X^i+
\bar\zeta_a(\dot\zeta^a+A^a_{bi}(X)\dot X^i\zeta^b)
\right)dt
$$
Here, $A$ represents the connection 1-form of $\nabla^{E}$. $\zeta$ and $\bar\zeta$ form a conjugate pair of Fermionic variables associated with $E$, so their kinetic term should naturally be $\bar\zeta\nabla^E\zeta$. The classical equations of motion are:

$$\dot X=0, \qquad\dot\zeta=\dot{\bar\zeta}=0. $$

\par Next, we extend the action to phase space, namely $T^*M$, where the momentum components $p_i$ serve as fiber coordinates. The Poisson brackets are defined as follows:
$$
\{X^i,p_j\}_{\mathrm P}=\delta^i_j,\qquad
\{\zeta^a,\bar\zeta_b\}_{\mathrm P}=\delta^a_b.
$$

\par We are already familiar with the method of quantizing this pair of fermions: after quantization, the annihilation operator $\zeta^a$ and creation operator $\bar\zeta_b$ act on their respective state spaces in a manner analogous to $\iota_{e_b}$ and $de_a\wedge$ acting on $\Lambda^*E|_U$. We consider only those operators for which, given a matrix $B$, we define $\mu(B)=\bar\zeta_b B^b_a \zeta^a$ or $(de_b\wedge) B_a^b 
(\iota_{e_a})$, depending on the context. Such operators commute with the particle number operator $N=\bar\zeta_a\zeta^a$; consequently, it is entirely permissible to restrict the entire state space $\Lambda^*E|_U$ to the one-particle state space $\Lambda^1E|_U\simeq E|_U$. Simple calculations show that the action of $\mu(B)$ on $\Lambda^1E|_U$ is equivalent to the action of $B\in \End(E|_U)$ on $E|_U$.
Calculations yield
$$
\{\mu(B),\mu(C)\}_{\mathrm P}=\mu([B,C]).
$$

\par The free action in phase space that we expect to observe is given by:
$$
S[X,p,\bar\zeta,\zeta]=-\int\left(p_i\dot X^i+\bar\zeta\dot\zeta\right)dt
$$
It is observed that as long as it satisfies constraint
$$
\theta_i=\lambda_i(X)-p_i+\mu(A_i(X))=0,
$$
the action we initially introduced is obtained.

\subsection{Handling constraints}

\par The aforementioned set of constraints satisfies
$$
\{\theta_i,\theta_j\}_{\mathrm P}
=\omega_{ij}+\hbar\mu(R^E_{ij}),
$$
where the curvature of $\nabla^E=d+A$ is given by $R^E=dA+A\wedge A$. Since this constitutes a set of second-class constraints, we introduce additional variables $y^i$ to reduce them to first-class constraints. These variables actually represent the fiber coordinates of $TM$, and they must satisfy
$$
\{y^i,y^j\}_{\mathrm P}=\pi^{ij}=(\omega^{-1})^{ij}.
$$

\par First, we seek a Hamiltonian function that generates the symplectic connection $\nabla^{TM}=d+\Gamma$:
$$
G_i=\frac12\omega_{ab}\Gamma^b_{ic}y^a y^c,\qquad G=G_idX^i.
$$
When the connection $\nabla^{TM}$ is induced onto the completed symmetric algebra of the cotangent bundle, denoted as $\widehat{Sym}(T^*M)=:W$, we obtain
$$
\nabla^W a = da + \{G,a\}_{\mathrm P}.
$$
Next, we determine the covariant variables, which correspond to the momentum generating the total covariant derivative:
$$
\Pi_i = p_i - G_i - \mu(A_i).
$$
Given that the full connection is $\nabla=\nabla^W\otimes\mathbb{1}+\mathbb{1}\otimes\nabla^E$, we have
$$
\nabla a = da + \{G+\mu(A),a\}_{\mathrm P},
\qquad
\{\Pi_i,a\}_{\mathrm P} = -\nabla_i a.
$$
In particular,
$$
\nabla\mu(B) = \mu(dB+[A,B]).
$$
Through further computation, it follows that
$$
(\nabla)^2 a = \{R+\mu(R^E),a\}_{\mathrm P},
$$
where
$$
R = dG + \frac12\{G,G\}_{\mathrm P} = \frac14\omega_{ab}\mathcal R^b{}_{cij}\,y^a y^c\,dX^i\wedge dX^j.
$$
$\mathcal R^b{}_{cij}$ is the curvature of $\nabla^{TM}$.

\par Next, we introduce correction terms involving $ y $ to modify $ \theta $ such that it becomes the first-class constraint. Firstly, it is natural to replace the momentum $p$ with the covariant momentum $\Pi$. Subsequently, by attempting this approach in the absence of curvature, it is found that modifying the aforementioned constraint to $T_i = \lambda_i - \Pi_i + \omega_{ji} y^j$ renders it of the first class. Therefore, we define
$$
\gamma = \gamma_idX^i,
\qquad
\gamma_i=\omega_{ji} y^j + r_i,
\qquad
T_i = \lambda_i - \Pi_i + \gamma_i.
$$
It can then be shown that the condition $\{T_i, T_j\} = 0$ for all $i, j$ is equivalent to
$$
R + \mu(R^E) + \nabla\gamma + \frac{1}{2}\{\gamma, \gamma\}_{\mathrm P} = -\omega,
$$
which represents the classical version of the Fedosov's equation.

\par Recall that at this stage, the system possesses the gauge symmetry generated by $\{T_i,\quad\}_P$.

\subsection{Deformation quantization}

\par Canonical quantization in Euclidean space:
$$[\hat x^i, \hat p_j]=i\hbar\delta^i_j,
\qquad
[\hat x^i, \hat x^j]=[\hat p_i, \hat p_j]=0.$$
Inspired by this approach, deformation quantization directly replaces the multiplication between functions with a star product, which corresponds to the composition of operators obtained from quantizing those functions. To be precise, this is not feasible, but the star product attempts to emulate the behavior of the formula above. In flat symplectic spaces, this operation is standard and is known as the Moyal product:
$$F\star G=m\exp\left(\frac{\hbar}{2}\pi^{ij}\partial_i\otimes\partial_j\right)F\otimes G$$
On symplectic manifolds, greater flexibility arises, with the general requirement being:
$$
F\star G=FG+\frac{\hbar}{2} \{F,G\}+\hbar^2 C_2(F,G)+\cdots.
$$

\par We first perform quantization in local coordinates within a local coordinate chart, that is, by introducing the star product. Its action on $y$ is analogous to the Moyal product:
$$a\star b=m\exp\left(\frac{\hbar}{2}\pi^{ij}\partial_{y^i}\otimes\partial_{y^j}\right)a\otimes b.$$
In terms of $\zeta$ and $\bar \zeta$, it adopts the form of canonical quantization:
$$[\bar\zeta_b,\zeta^a]_\star=\hbar\delta^a_b, \qquad[\bar\zeta_b,\bar\zeta^a]_\star=[\zeta_b,\zeta^a]_\star=0.$$
At this point,
$$
\frac1\hbar[\mu(B),\mu(C)]_\star=\mu([B,C]).
$$
Its action on $X$ and $p$ also corresponds to their canonical quantization:
$$
[p_i,a(X)]_\star=-\hbar\partial_{X^i} a(X),\qquad[p_i,p_j]_\star=0.
$$
Or more precisely, formulate as:
$$
\bar\zeta_a\star\zeta^b=\bar\zeta_a\zeta^b,
\qquad
\zeta^b\star\bar\zeta_a
=-\bar\zeta_a\zeta^b+\hbar\delta^b_a
$$
$$
a\star b=m\exp\left(\hbar\partial_{X^i}\otimes\partial_{p_i}\right)a\otimes b
$$
Now the action should be understood as $\bar\zeta_a=d e_a\wedge, \zeta^b=\hbar\iota_{\partial_{e_b}}$.

\par It can be observed that this approach is equivalent to adopting $X,p$ as Darboux coordinates, and performing a Moyal product when the symplectic form is given by $dX^i\wedge dp_i$. Thus, the motivation behind the Moyal product or deformation quantization lies in attempting to emulate canonical quantization. However, this construction can also be derived in a more reasonable manner. Consider a one-dimensional field theory:
$$\R\to (\mathbb{R}^{2n},\omega)$$
Treating $(\mathbb{R}^{2n},\omega)$ as the phase space, the free action is
$$S[y(t)]=\frac{1}{2}\int_{\R} y^i\omega_{ij}\dot y^j dt$$
which is precisely equivalent to $S=\int pdq$. For an observable $\mathcal O (y)$, its expectation value can be obtained via path integration (normalization issues are considered separately for now):
$$\left<\mathcal O\right>=\int D[y(t)]\mathcal O (y)e^{\frac{1}{\hbar}S[y]}$$
First, computing the correlation function of coordinate functions yields:
$$\left<y^i(t_1)y^j(t_2)\right>=\hbar\pi^{ij}G(t_1,t_2)$$
where
$G(t_1,t_2)=\frac12\operatorname{sgn}(t_1-t_2)$
represents the Green’s function of $\frac{d}{dt}$ on $\R$. Further calculations show that for observables $f,g$ which we regard as formal polynomials in the coordinates, given a fixed point $Y$, we have:
$$
\begin{aligned}
  \left<f(Y+y(t_2))g(Y+y(t_1))\right>
  &=m\exp\left(\frac{\hbar}{2}\pi^{ij}\partial_{Y^i}\otimes\partial_{Y^j}\right)f(Y)\otimes g(Y) \\
  &=f\star g(Y).
\end{aligned}
$$
The computation proceeds by Taylor expanding $f(Y+y(t))$ at $y(t)=0$ with respect to the variable $y(t)$, and then applying the Wick contractions of the coordinate functions given above. Thus, in a sense, what we are computing here is the effect of quantum fluctuations around a fixed point. In other words, when this free phase-space action is quantized through path integration, The Moyal product can be interpreted as the effective observable associated with the joint measurement of two observables. Henceforth, we shall always place operators evaluated at later times to the left.

\par The “quantization” process itself involves merely two steps: first, select any closed 2-form with magnitude $\hbar$ to serve as the formal symplectic form:
$$\Omega_\hbar=\omega+\hbar\omega_1+\cdots=d\lambda_\hbar$$
Second, replace $\{\quad,\quad\}_P$ with $[\quad,\quad]_\star$. Consequently, the constraint condition $T$ satisfies $\frac{1}{\hbar}[T_i,T_j]_\star=0$ after quantization, which is equivalent to
$$
R + \mu(R^E) + \nabla\gamma + \frac{1}{2\hbar}[\gamma, \gamma]_\star = -\Omega_\hbar
$$
This equation constitutes the Fedosov equation. Fedosov showed that for any closed 2-form $\Omega_\hbar$, a unique $\gamma$ can be obtained via recursive computation of $r$. In this context, $\gamma$ is a 1-form on a manifold dependent on $\hbar$ and taking $X$, $y$, $\zeta$, and $\bar\zeta$ as its variables, where $\zeta$ and $\bar\zeta$ always appear in the form of certain $\hbar\mu(B)$.

\subsection{Quantum BRST algebra and gauge symmetry}

\par Next, the ghosts $C^i, P_i$ are introduced to constrain $T_i$, satisfying
$$[C^i, P_j]_\star = \hbar\delta^i_j,\qquad [P_i, P_j]_\star = [C^i, C^j]_\star = 0.$$
The BRST charge $Q = C^i T_i$ generates the quantum BRST differential $s = \frac{1}{\hbar}[Q, \cdot]_\star$. The condition $[Q, Q]_\star = 0$ ensures that $s^2 = 0$. We have:
$$
s C^i = 0,\qquad s P_i = T_i,\qquad s T_i = 0.
$$
Further calculations show that the action of $s$ on $X, y$ is given by
$$
s a = C^i\left(\nabla_i a + \frac{1}{\hbar}[\gamma, a]_\star\right).
$$
In other words, if $C^i$ is identified with $dX^i$, then $s$ corresponds to the Fedosov connection $D = \nabla + \frac{1}{\hbar}[\gamma, \cdot]_\star$. This reveals the intrinsic meaning of the Fedosov connection: since $\frac{1}{\hbar}[T_i, \cdot]_\star$ generates gauge symmetry within this quantum system, and $dX^i$ serves as the gauge parameter, $D$ functions as the gauge transformation operator.

\par At this stage, the variables of the quantum BRST algebra—i.e., the functions in the system—are $X, y, p, C, P, \zeta, \bar\zeta$, which constitute the space in which quantum observables exist. It is known that $p$ and $y$ are respectively the coordinates on the fibers of $T^*M$ and $TM$. As discussed earlier, since the state space after quantizing fermionic variables is restricted to a single particle, fermionic operators appear only in the form of certain $\mu(B)$; accordingly, $\zeta$ and $\bar\zeta$ appear only as certain $\hbar\mu(B)$. If $C^i$ is identified with $dX^i$, then the corresponding $P_i$ must be some odd derivative $\partial_i$. Thus, the quantum BRST algebra can be identified as
$$
\Omega_M^*\left(\widehat{Sym}(TM) \otimes \End(E) \otimes \widehat{Sym}(T^*M) \otimes \Lambda^*(TM)\right)[[\hbar]].
$$
This implies that we need to establish, or verify, the transformation rules governing these variables with respect to coordinates. According to general BRST theory, in the zeroth cohomology of the quantum BRST algebra with respect to $s$, both $P$ and $T$ vanish; $P_i\partial_{T_i}$ serves as the homotopy operator between them. Subsequently, after imposing the constraint conditions, the phase space variable $p$ also disappears. Hence, the zeroth cohomology space of the quantum BRST algebra with respect to $s$—i.e., the set of quantum observables satisfying the constraints and gauge conditions—is precisely the $D$-flat sections of Weyl bundle $W \otimes \End(E)[[\hbar]]$.

\par For simplicity, we shall first provide a brief explanation using the Euclidean space prior to quantization. In this case, $r = 0$, and the gauge transformation takes the form $\delta = \varepsilon^i\{T_i,\quad\}$. The calculations yield
$$
\delta X^i = \varepsilon^i, \qquad \delta y^i = -\varepsilon^i.
$$
Thus, $X + y$ serves as a gauge-invariant quantity under the gauge symmetry we have been studying. Consequently, the set $\{y = 0\}$ intersects every gauge orbit exactly once. This implies that each gauge equivalence class contains exactly one representative element independent of $y$. The result proved by Fedosov,
$$
H^0(\Gamma(M,W \otimes \End(E))[[\hbar]], D) \simeq \Gamma(M,\End(E))[[\hbar]],
$$
reflects precisely this fact. Therefore, we first embed the functions $f$ and $g$ on the manifold into the quantum BRST algebra and seek genuine physical observables—elements closed under gauge transformations—$F$ and $G$ within their respective gauge equivalence classes. Since we have already quantized the BRST algebra by introducing the star product, we then compute this star product $F \star G$ between the physical observables. Finally, we select an element of $y = 0$ from the gauge equivalence class of $F \star G$ and define it as $f \star g$. This is precisely deformation quantization. Here, we introduce a pair of Fermionic variables, thereby obtaining directly the case with vector bundles.

\section{Algebraic index theorem in BV formalism}
 
\subsection{Extended action}

\par We first return to the stage prior to deformation quantization, where $T_i = \lambda_i - \Pi_i + \gamma_i, s=C^i\{T_i,\quad\}_P,\{T,T\}_P=0$. Following the BV formalism, we incorporate the constraints into the action by means of Lagrange multipliers, introduce antifields for all field variables $\Phi^A=(X,p,y,\zeta,\bar\zeta,C,u)$, and further introduce the auxiliary fields $\bar C,B$ together with their antifields. This yields the extended action
$$S_{\mathrm{BV}}
=\int_0^1\left[
-p_i\dot X^i
-\frac12\omega_{ij}y^i\dot y^j
-\bar\zeta_a\dot\zeta^a
+u^iT_i+\Phi^*_As\Phi^A+\bar C^*_i B^i
\right]dt$$
which can, of course, be written out completely explicitly. Here we have additionally included the kinetic term $\frac12\omega_{ij}y^i\dot y^j$. This is a natural choice, since the preceding discussion suggests that $y$ can be regarded as a certain expansion or extension of $X$, and should therefore possess a kinetic term corresponding to the $p\dot X$ term. To avoid confusion with the notation used in subsequent sections, we shall denote the previously defined BV bracket by $(\quad,\quad)_{BV}$. One can verify that this action satisfies the classical master equation:
$$\frac12(S,S)_{BV}=0$$
or equivalently, $s_{BV}^2=0$. It also formally satisfies the quantum master equation
$$\frac12(S,S)_{BV}+\hbar \Delta S=0$$
since $\Delta S=0$.

\subsection{Gauge fixing} 

\par We now perform gauge fixing. We introduce the background variables $x^i,c^i$. These can be understood as the time-independent counterparts of $X(t),C(t)$ on the manifold, corresponding respectively to a point $x$ on the manifold and the differential forms $dx$. We choose the gauge-fixing fermion
$$
\Psi=\int_0^1\bar C_i(t)(X^i(t)-x^i)\,dt.
$$
The corresponding Lagrangian submanifold $\mathcal L_\Psi$ is specified by
$$
\begin{gathered}
X_i^*=-\bar C_i,
\qquad \bar C^{*i}=-(X^i-x^i),
\\
p^*=u^*=C^*=B^*=0.
\end{gathered}
$$
A direct calculation shows that restricting the extended action to this Lagrangian submanifold is equivalent to adding the following two terms to the original antifield-independent action:
$$(x_i^*-\int_0^1\bar C_i\,dt,)C^i
$$
$$
s\Psi_X=\int_0^1\left[B_i(X^i-x^i)-\bar C_i(C^i-c^i)\right]dt.
$$
We subsequently denote the first term simply by $x_i^*C^i$. In the path integral, integration over $B$ and $\bar C$ then imposes the conditions
$$X(t)=x,\qquad C(t)=c.$$

\par The gauge-fixing fermion introduced above does not involve the variables $\zeta,\bar \zeta,y$. In particular, $y$ is an auxiliary variable introduced to resolve the constraints and does not appear in the standard BV formalism; consequently, the preceding gauge-fixing procedure does not involve it. We therefore need to choose the additional Lagrangian submanifold $\mathcal L=\{\zeta^*=\bar \zeta^*=0\}$ explicitly. For convenience, we introduce the notation
$$y_a^*=\omega_{ab}\eta^b,
\qquad
d_y=\eta^a\partial_{y^a}.$$
A lengthy direct calculation then yields
$$y^*_asy^a=-c^id_y(G_i+\mu(A_i)+\gamma_i).$$
Here $d_y\mu(A_i)=0$. The action thus reduces to
$$S_{gf}=x^*_ic^i+\int_0^1\left[-\frac12\omega_{ab}y^a\dot y^b
-\bar\zeta_\alpha\dot\zeta^\alpha-c^i\eta^a\frac{\delta}{\delta y^a}\left(G_i+\gamma_i\right)\right]dt. $$

\par We now need to incorporate quantum corrections. As demonstrated in the preceding chapter, the free part of this action already gives rise to quantum contractions governed by the star product rather than ordinary multiplication. We must therefore replace the classical constraint condition $\{T,T\}_P=0$ with its quantum counterpart $[T,T]_\star=0$. In the present setting, this amounts to replacing $\gamma$, obtained by solving the classical Fedosov equation associated with $\{T,T\}_P=0$, with the solution $\gamma$ of the genuine Fedosov equation corresponding to $[T,T]_\star=0$. Equivalently, the replacement of $\gamma$ can be understood as compensating for the additional quantum effects introduced by replacing the classical constraint condition $\{\quad,\quad\}_P=0$ with its quantum counterpart $[\quad,\quad]_\star=0$. From this perspective, it may also be regarded as part of the quantum correction to the action required by BV quantization. However, this modification has the unintended consequence that the resulting action no longer satisfies the quantum master equation (QME). Consequently, integration over a chosen Lagrangian submanifold is no longer necessarily justified within the BV formalism. Nevertheless, we shall proceed with the path integral computation along these lines. The calculations presented below should therefore be regarded, to some extent, as formal and heuristic physical arguments, even though they ultimately reproduce the exact result. The action thus requires further quantum corrections. In fact, the correct construction is already within reach: after performing the integrations described in the following two sections, the resulting effective action, which depends only on the remaining background fields, satisfies the quantum master equation. Along this way, \cite{GLL17} employs the homotopy renormalization group flow to construct, starting from an action of the type considered here, a family of modified effective actions satisfying the QME. The path integral is then replaced by taking some limit of the parameter of this family, thereby rendering the procedure mathematically rigorous.This physically and mathematically consistent construction yields the same result as our naive path integral calculation.

\par In accordance with our previous discussion of the representation of $\zeta,\bar\zeta$ after quantization, we make the change of variables
$$
\bar\zeta_a=\bar z_a,
\qquad \zeta^a=\hbar z^a.
$$
Correspondingly, $\gamma(x,y,\zeta,\bar\zeta)$ should be replaced by $\frac1\hbar\gamma(x,y,z,\bar z)$.

\par Next, we decompose $y,\eta$ into their time-independent zero modes and the remaining fluctuations:
$$
y(t)=y+\widetilde y(t),\qquad
\eta(t)=\eta+\widetilde\eta(t)
$$
where $\int_0^1\widetilde y dt=
\int_0^1\widetilde\eta dt=0$. We choose the Lagrangian submanifold
$$
\mathcal L
=\{\widetilde\eta=0\}
$$
At this stage, we have specified the Lagrangian submanifolds for all field variables, and the final path integral is performed over their intersection, a real Lagrangian submanifold in whole space. Strictly speaking, the treatment of the background-field is not yet complete; we shall return to this issue later. In particular, since $c^id_yG_i$ is linear in $y$, its integral reduces to
$$\int_0^1c^id_yG_i(x,y+\widetilde{y}(t),\eta)dt=c^id_yG_i(x,y,\eta). $$

\par Recall that, since the action of $z,\bar z$ is restricted to the one-particle sector, the dependence of $\gamma$ on $z,\bar z$ occurs exclusively through bilinear terms of the form $\bar z_aB^a_bz^b$ (which are moreover of at least first order in $\hbar$). If we regard $\gamma$ as a solution of the Fedosov equation, the corresponding term is simply $B\in\End(E)$. We shall adopt this identification in what follows. Thus, the gauge-fixed action can now be written as
$$
S[x,x^*,c,y,\widetilde{y}(t),\eta,z(t),\bar z(t)]=x^*_ic^i-c^id_yG_i+\int_0^1\left[-\frac{1}{2}\omega_{ij}\widetilde{y}^i\dot{ \widetilde{y}^j}-\hbar \bar z(\partial_t+\frac{c^id_y\gamma_i}{\hbar})z\right]dt. 
$$
Note that $\gamma_i$ still depends on the full variable $y+\widetilde{y}$.

\par At this stage, we can make the following geometric identification. Consider the so-called local quantum observables: such observables do not depend on the entire field $y(t)$, but only formally on $y$ and a fixed time $t_0$, as if an observable were inserted locally at a single point on the worldline. These observables depend on $x,c,y,\eta,\bar z Bz$, and therefore belong to
$$\Omega^*_M\otimes\widehat{Sym}(T^*M)\otimes \Lambda^{-*}(T^*M)\otimes\End(E)$$
where the negative grading on $\Lambda^{-*}(T^*M)$ reflects the fact that $\eta$ -- we identify it with $dy$ -- has ghost number $-1$. We shall, of course, restrict our attention to the original observables that are independent of the antifield $\eta$. Although we have restricted the extended action to a Lagrangian submanifold, the action of $\{S,\quad\}_{BV}$ must be computed using the full action before restricting the result to the submanifold, just as in ordinary differential geometry. This makes the calculation rather lengthy. Roughly speaking, the terms $x^*_ic^i,-c^id_yG_i,\bar\zeta^*A\zeta,kinetic\,terms$ give rise to $d_M,\Gamma,[A,\quad],\partial_t$, respectively. The quantum differential $\{S,\quad\}_{BV}+\hbar\Delta$ then acts on local observables as $\nabla+\{\gamma,\quad\}_P$. We subsequently incorporate the quantum correction by replacing this operator with $\nabla+\frac1\hbar[\gamma,\quad]_\star$, which is precisely the Fedosov connection. Consequently, the space of local observables satisfying the gauge-invariance condition is identified with the space of flat sections of the Weyl bundle with coefficients in $\End(E)$ with respect to the Fedosov connection. In what follows, we shall consider only observables of this type.

\subsection{Path integral of internal term}

\par We first compute the path integral over the internal degrees of freedom associated with $E$:
\[
\int
D\bar z\, Dz\;
e^{-\int_0^1\bar z(\partial_t+\frac{c^id_y\gamma_i}{\hbar})z\,dt}
\]
In what follows, we abbreviate $\frac{c^id_y\gamma_i}{\hbar}$ as $B$. We restrict our attention to field configurations defined on loops, so that $\int_0^1$ is replaced by $\int_{S^1}$. The resulting path integral therefore computes the partition function. Here we adopt the method of regularized functional determinants commonly used in supersymmetric quantum mechanics.

\par In the present setting, $B$ is time-dependent because $\gamma$ depends on $\widetilde{y}$. Consider the differential equation
\[
\dot U(t)=-B(t)U(t),\qquad U(0)=I.
\]
Its solution is given by
\[
U(t)=T\exp\!\left(-\int_0^tB(s)\,ds\right).
\]
Here $T$ denotes the time-ordering operator, which prescribes a rearrangement of the products arising from the expansion of the exponential, with operators evaluated at later times placed to the left. If we replace $B$ by its corresponding operator $\mu(B)$, then $U(t)$ is precisely the time-evolution operator of this system. Let $K=-\log U(1)$ and $g(t)=U(t)e^{tK}$. The latter defines a transformation satisfying
\[
g^{-1}(t)\left(\partial_t+B(t)\right)g(t)=\partial_t+K.
\]
Thus, the original time-dependent operator is transformed into one with constant coefficients.

\par We first compute the determinant ratio $\det(\partial_t+K)/\det(\partial_t)$. In the present setting, since the trace $\operatorname{tr}$ on $\End(E)$ is the ordinary trace rather than the supertrace, we impose antiperiodic boundary conditions on $z,\bar z$ over $S^1$. The corresponding antiperiodic Fourier modes are
\[
e^{2\pi i(m+\frac12)t},\qquad m\in\mathbb Z.
\]
A pair of modes with frequencies $\pm(2m+1)\pi$ then contributes
\[
\frac{
\det_E(i(2m+1)\pi I+K)
\det_E(-i(2m+1)\pi I+K)
}{
\det_E(i(2m+1)\pi I)
\det_E(-i(2m+1)\pi I)
}
=
\det\left(I+\frac{K^2}{(2m+1)^2\pi^2}\right).
\]
Consequently, the functional determinant is given by
\[
\begin{aligned}
\frac{
\operatorname{Det}(\partial_t+K)
}{
\operatorname{Det}(\partial_t)
}
&=
\prod_{m=0}^{\infty}
\det\left(
I+\frac{K^2}{(2m+1)^2\pi^2}
\right) \\
&=\det\left(\cosh\frac{K}{2}\right).
\end{aligned}
\]

\par We next apply zeta-function regularization. We first consider the one-dimensional operator $\partial_t$. For
\[
\operatorname{Det}(\partial_t)
=
\prod_{n=0}^{\infty}
\pi^2(2m+1)^2.
\]
consider the function
\[
\begin{aligned}
  F(s)
=
\pi^{-2s}\sum_{n=0}^{\infty}(2n+1)^{-2s}
=
\pi^{-2s}(1-2^{-2s})\zeta(2s),
\end{aligned}
\]
Using the special value $\zeta(0)=-\frac{1}{2}$ of the Riemann zeta function, we obtain the regularized product
\[
\prod_{n=0}^{\infty}
\pi^2(2n+1)^2
=\exp(-F'(0))=2.
\]
Thus, for a vector bundle of rank $r$, we have $\operatorname{Det}(\partial_t)=2^r$. Combining this with the identity
$
2\cosh\frac x2=e^{x/2}(I+e^{-x})
$, we obtain
\[
\operatorname{Det}(\partial_t+B)
=
e^{\frac12\operatorname{tr}_E K}
\det_E(I+U(1)).
\]
Finally, since
\[
\frac d{dt}\log\det_E U(t)
=
\operatorname{tr}_E\bigl(U^{-1}\dot U\bigr)
=
-\operatorname{tr}_E B(t),
\]
we conclude that
\[
\operatorname{tr}_E K
=
\int_0^1\operatorname{tr}_E B(t)\,dt.
\]

\par We now need to make two further adjustments. The preceding calculation was performed for a general pair of complex fermionic variables, whereas we wish to restrict their representation to the one-particle sector. First, in quantizing $\zeta,\bar\zeta$, we adopted the normal-ordering prescription $\mu(B)=\bar\zeta_aB^a_b\zeta^b$, while the functional determinant calculation above, without any additional prescription, corresponds to symmetric ordering $\mu(B)=\frac12(\bar\zeta_aB^a_b\zeta^b-\zeta^bB^a_b\bar\zeta_a)$. These two prescriptions differ by $\frac12\tr(B)I $. We must therefore remove the additional factor
$
\exp\!\left(\frac12\int_0^1\operatorname{tr}_E B(t)\,dt\right)
$
from the determinant obtained above.

The second adjustment is to restrict the result to the one-particle representation. To this end, we introduce a formal parameter $\nu$, so that $\mu(B-\nu I)=\mu(B)-\nu N$, where $N$ is the particle-number operator. The evolution operator $U_{B-\nu I}$ associated with $\partial_t+B-\nu I$ is related to the evolution operator $U_B$ associated with $\partial_t+B$ by
$$U_{B-\nu I}(t)=e^{\nu t N}U_B(t).$$
Thus, introducing $\nu I$ serves only to keep track of the particle number, without otherwise modifying the evolution. We therefore consider $\partial_t+B-\nu I$ instead. Extracting the term of first order in $e^\nu$ from
$$\det(I+U_{B-\nu I}(1))=\det(I+e^\nu U(1))=1+e^\nu \tr_E U(1)+O((e^\nu)^2)$$
then isolates precisely the contribution from the one-particle sector.

\par Finally, recall that for a pair of complex fermionic variables $z,\bar z$, the formal Gaussian integral satisfies
$$\int Dz D\bar z e^{-\bar z A z}=\det(A).$$
Combining all the preceding arguments in this section, we obtain the path integral over the internal degrees of freedom restricted to the one-particle representation:
\[
\int
D\bar z\, Dz\;
e^{-\int_{S^1}^{N=1}\bar z(\partial_t+\frac{c^id_y\gamma_i}{\hbar})z\,dt}=\tr_E T e^{-\int_{S^1}d_y\gamma/\hbar\, dt}. 
\]

\subsection{Path integral of non-zero modes}

\par We next evaluate the path integral over $\widetilde{y}$:
$$\int D[\widetilde{y}] \left(\tr_E T e^{-\int_{S^1}d_y\gamma/\hbar\, dt}\right)e^{\frac1\hbar\int_{S^1}-\frac{1}{2}\omega_{ij}\widetilde{y}^i\dot{ \widetilde{y}^j}}$$
The computation proceeds essentially as in the preceding chapter. Locally, consider the field theory
$$S^1 \to (\mathbb{R}^{2n},\omega)$$
with free action $\int_{S^1}-\frac{1}{2}\omega_{ij}y^i\dot{ y^j}$. The corresponding two-point correlation function is
$$\left<y^i(t_1)y^j(t_2)\right>=\hbar\pi^{ij}P(t_1,t_2)$$
where
\[
P(t_1,t_2)=t_2-t_1+\frac12,
\qquad
P(t,t)=0.
\]
Here $P$ is the Green's function for $\frac{d}{dt}$ on $S^1$. For two local observables $f,g$, we then have
$$
\begin{aligned}
  \lim_{t_1\to t_2}\left<f(Y+y(t_2))g(Y+y(t_1))\right>
  &=m\exp\left(\frac{\hbar}{2}\pi^{ij}\partial_{Y^i}\otimes\partial_{Y^j}\right)f(Y)\otimes g(Y) \\
  &=f\star g(Y).
\end{aligned}
$$
This admits a natural local interpretation. When two local observables collide on $S^1$, a sufficiently small neighborhood $U\subset S^1$ of the collision point is locally indistinguishable from $\R$. Consequently, the short-distance behavior of their correlation function reproduces the Moyal product obtained from the field theory on $\R$ discussed earlier. Thus, in the present model, the Moyal product can be interpreted as the effective observable resulting from the collision and fusion of two local observables into a single one. Note that by ``local observables'' we mean observables inserted at individual points on the worldline, rather than observables that depend on the entire worldline in the conventional sense. Consequently, their expectation values still depend explicitly on the times at which they are inserted. We choose the counterclockwise orientation as the positive orientation of $S^1$. For a collection of observables $f_0(t_0),f_1(t_1),\cdots,f_m(t_m)$ arranged in clockwise order along the circle, Wick's theorem gives their joint expectation value as
$$\left\langle f_0(Y+y(t_0))\cdots f_m(Y+y(t_m))\right\rangle=\exp\left(\frac{\hbar}{2}\pi^{ij}P(t_a,t_b)\partial_{Y^i_a} \partial_{Y^j_b}\right)f_0(Y_{0})\cdots f_m(Y_{m})\big|_{Y_{a}=Y}$$
Here the indices $\{a,b,\dots\}$ range over $\{0,1,\dots,m\}$.

\par We adopt the notation $\langle\quad\rangle$ to denote expectation values with respect to the free action $\int_{S^1}-\frac{1}{2}\omega_{ij}\widetilde{y}^i\dot{ \widetilde{y}^j}$ for the field $\widetilde{y}(t)$. Thus, the quantity we now need to compute is
$$\left\langle\tr_E T e^{-\int_{S^1}d_y\gamma/\hbar\, dt}\right\rangle.$$
Note that the decomposition $Y+y(t)$ in the preceding example is precisely the decomposition $y+\widetilde{y}(t)$ of the field $y(t)$ appearing in $\gamma$ in the present setting. We reconsider the evolution operator introduced earlier:
\[
U(t)=T\exp\!\left(-\int_0^tB(s)\,ds\right).
\]
By iterating the integral equation satisfied by $U(t)$, we obtain its explicit expansion, known as the Dyson series:
\[
\begin{aligned}
U(1)
={}&\mathbb{1}-\int_0^1B(t_1)\,dt_1\\
&+\int_0^1dt_2\int_0^{t_2}dt_1\,
B(t_2)B(t_1)\\
&-\int_0^1dt_3\int_0^{t_3}dt_2
\int_0^{t_2}dt_1\,
B(t_3)B(t_2)B(t_1)+\cdots .
\end{aligned}
\]
Equivalently,
\[
T e^{-\int_0^1B(t)\,dt}
=
\mathbb{1}+
\sum_{k\geq 1}(-1)^k
\int_{\Delta_k}
B(t_k)\cdots B(t_1)\,d\mathbf t.
\]
Here $\Delta_k$ denotes the $k$-dimensional simplex
\[
\Delta_k
=
\{(t_1,\ldots,t_k)|0<t_1<\cdots<t_k<1\}.
\]

\par We therefore obtain
$$
\begin{aligned}
  &\left\langle\tr_E T e^{-\int_{S^1}d_y\gamma/\hbar\, dt}\right\rangle \\
  =&\Sigma_{k\geq 0}(-1)^k\int_{S^1\times \Delta_k}dt_0dt_1\cdots dt_k\,\tr_E\,m\,e^{\frac{\hbar}{2}\pi^{ij}P(t_a,t_b)\partial_{y^i_a}\otimes\partial_{y^j_b}}\left(\mathbb{1}\otimes (d_y\gamma/\hbar)^{\otimes k}\right)|_{y_a=y}
\end{aligned}
$$
Here we have introduced an additional insertion of $\mathbb{1}$ at a point on $S^1$ to provide a designated insertion point for observables in subsequent calculations. Equivalently, what we are computing here is the expectation value of the identity observable, namely the partition function. The expression above is usually understood as arising from the construction of a local expectation functional, or equivalently, an expectation functional at the Lie algebra level, followed by its descent to a global expectation functional using the connection one-form $\gamma$, in accordance with the philosophy of formal geometry. (ref.\cite{NT95}\cite{GLX21}). Our calculation, however, demonstrates that this construction emerges naturally and directly from the physical path integral. 

\par We denote the resulting expression by $e^{\gamma_\infty/\hbar}$. In other words, we interpret $\gamma_\infty$ as the effective action obtained by integrating out the fields $\widetilde{y},z,\bar z$ in the gauge-fixed action, while retaining the time-independent background terms. For the total action $S=x^*c-cd_y G+\gamma_\infty$ at this stage, one can verify that it satisfiesthe quantum master equation $\frac12\{S,S\}_{BV}+\hbar\Delta S=0$. 

\subsection{Integration over the background fields}

\par The action now takes depends only on the time-independent background fields, or equivalently, the zero modes $x,x^*,c,y,\eta$. As mentioned earlier, one issue remains to be addressed: the gauge-fixing Lagrangian submanifolds previously chosen for $x,x^*,y,\eta$ involve only their nonzero modes, namely the time-dependent fluctuations. The time-independent background components, or zero modes, have not yet been gauge-fixed. Consequently, the zero-mode antifields $x^*,\eta$ are still present. We now choose the Lagrangian submanifold
$$\mathcal L=\{x^*=0,\quad y=0\}.$$
Only at this stage have we genuinely restricted the extended action to a Lagrangian submanifold of the full space of fields. The action then reduces to $\gamma_\infty|_{y=0}$. Integration over $c$ is likewise a Berezin integral. Upon identifying $c^i$ with $dx^i$, its combination with integration over $x$ gives precisely the geometric integral $\int_M$. Integration over $\eta$ is a Berezin integral. Considering numbers of $c^i$, identifying $\eta$ with $dy$ and adopting the volume-normalization convention used in other literature, this integration can be represented by the operator $\exp(\hbar\frac12\pi^{ij}\iota_{y^i}\iota_{y^j})=e^{\hbar\iota_\pi}$ before $\int_M$.

\par We thus arrive at the partition function of the one-dimensional field theory $S=-\int p\dot x+\bar\zeta\dot\zeta$, describing a free particle on the symplectic manifold $(M,\omega)$:
$$\mathcal Z_E=\int_{M} e^{\hbar\iota_\pi}e^{\gamma_\infty/\hbar}|_{y=0}.$$
This expression agrees exactly with the formula obtained from the construction in, for instance, \cite{GLX21}. In other words, the partition function of this field theory is precisely the algebraic index.

%\appendix

%\section{none}

\bibliographystyle{amsalpha}
\bibliography{references}

\providecommand{\bysame}{\leavevmode\hbox to3em{\hrulefill}\thinspace}
\providecommand{\MR}{\relax\ifhmode\unskip\space\fi MR }
\providecommand{\MRhref}[2]{%
  \href{http://www.ams.org/mathscinet-getitem?mr=#1}{#2}
}
\providecommand{\href}[2]{#2}
\begin{thebibliography}{GLX21}

\bibitem[Fed96]{F96}
B.~Fedosov, \emph{{Deformation quantization and index theory}}, Akademie
  Verlag, 1996.

\bibitem[FFS05]{FFS05}
Boris Feigin, Giovanni Felder, and Boris Shoikhet, \emph{Hochschild cohomology
  of the weyl algebra and traces in deformation quantization}, Duke
  Mathematical Journal \textbf{127} (2005), 487--517.

\bibitem[GL01]{GL01}
M.~Grigoriev and S.~Lyakhovich, \emph{Fedosov deformation quantization as a
  brst theory}, Communications in Mathematical Physics \textbf{218} (2001),
  437--457.

\bibitem[GLL17]{GLL17}
R.~Grady, Q.~Li, and S.~Li, \emph{Batalin–vilkovisky quantization and the
  algebraic index}, Advances in Mathematics \textbf{317} (2017), 575--639.

\bibitem[GLX21]{GLX21}
Z.~Gui, S.~Li, and K.~Xu, \emph{Geometry of localized effective theories, exact
  semi-classical approximation and the algebraic index}, Communications in
  Mathematical Physics \textbf{382} (2021), 441--483.

\bibitem[NT95]{NT95}
R.~Nest and B.~Tsygan, \emph{Algebraic index theorem}, Communications in
  Mathematical Physics \textbf{172} (1995), 223–262.

\end{thebibliography}

\end{document}